\magnification 1200
\baselineskip .85truecm
\input amssym.tex

\def \zgran{\displaystyle}

\def \za{\alpha}
\def \zb{\beta}

\def \zd{\delta}

\def \zl{\lambda}
\def \zm{\mu}

\def \zp{\pi}

\def \zf{\varphi}


\def \zlma{\ell}

\def \zsu{\sum}

\def \zmm{\pm}

\def \zpu{\cdot}
\def \zpor{\times}

\def \zmei{\leq}
\def \zmai{\geq}
\def \zco{\subset}

\def \zpe{\in}



\def \zinf{\infty}

\def \zfl{\rightarrow}

\def \zbv{\mid}
\def \zdbv{\parallel}
\def \z/{\over}

\centerline{\bf Polynomials Maps and Even Dimensional Spheres}
\medskip

\centerline{ Francisco-Javier TURIEL}

\centerline {Geometr{\'\i}a y Topolog{\'\i}a, Facultad de Ciencias
Ap. 59,}

\centerline{ 29080 M{\'a}laga, Spain}

\centerline{ email: turiel@agt.cie.uma.es}
\bigskip

{\bf Abstract.}  We construct, for every even dimensional sphere $S^n$, $n\zmai 2$, and every odd integer $k$, a
homogeneous polynomial map $f: S^{n}\zfl S^{n}$ of Brouwer degree $k$ and algebraic degree $2\zbv k\zbv -1$.

\centerline{--------------------------------------}
\bigskip \bigskip

A {\it polynomial map} from $X\zco {\Bbb R}^{m}$ to $Y\zco {\Bbb R}^{r}$ is the restriction to $X$ of a polynomial map
$F: {\Bbb R}^{m}\zfl {\Bbb R}^{r}$ such that $F(X)\zco Y$. When each component of $F$ is homogeneous of degree $k$,
we will say that the polynomial map from $X\zco {\Bbb R}^{m}$ to $Y\zco {\Bbb R}^{r}$ is {\it homogeneous} of degree
$k$.  As usual  $S^n$ is the sphere on ${\Bbb R}^{n+1}$ defined by the equation $x_{1}^{2}+...+x_{n+1}^{2}=1$, in
short $\zbv\zbv x\zbv\zbv^{2} =1$, whereas $S_{r}^{n}$, $1\zmei r\zmei n$,
will be the differentiable manifold, diffeomorphic to $S^n$ , defined by the equation
$(x_{1}^{2}+...+x_{r}^{2})^{2}+ x_{r+1}^{2}+...+x_{n+1}^{2}=1$. In this work we show:
\bigskip

{\bf Theorem 1.} {\it Suppose $n$ even and $\zmai 2$. Let $k$ be an integer. Then:

\noindent (a) If $k$ is odd, there exists a homogeneous polynomial map from $S^n$ to $S^n$ of Brouwer degree $k$
and algebraic degree $2\zbv k\zbv -1$.

\noindent (b) If $k$ is even there exists, for each $2\zmei 2r\zmei n$, a polynomial map from
$S^{n}_{2r}$ to $S^n$ of Brouwer degree $k$.}
\bigskip

Representing elements of $\zp_{n}(S^{n})$ by polynomial maps is an old question [1] which was affirmatively solved
by Wood, in 1968, provided that $n$ is odd (theorem 1 of [3], see[4] as well for the complex sphere).
Nevertheless, as far as I know, this problem is still open for $n$ even; our theorem settles it when the Brouwer
degree is odd. In both cases the polynomial maps constructed are homogeneous; therefore the problem of representing
elements of $\zp_{n}(S^{n})$ by homogeneous polynomial maps is solved now, since only zero and the odd topological
degrees may be represented in this way  when $n$ is even [2].

The proof of theorem 1 of [3] makes use of a natural polynomial map of topological degree 2 (lemmas 4 and 5).
Nothing similar is known for $n$ even; however the polynomial map
$x\zpe {\Bbb R}^{3}\zfl (x_{1}^{2}-x_{2}^{2}, 2x_{1}x_{2}, x_{3})\zpe {\Bbb R}^{3}$ send $S_{2}^{2}$ into $S^2$
with topological degree 2. Part (b) of our theorem generalizes this map.

For proving theorem 1 we start constructing a family of real polynomials in one variable. Let
$\zgran\zf_{\zlma}=\zsu_{j=0}^{\zlma}a_{j}t^{j}$ be the Taylor expansion of $\zf=(1-t)^{-1/2}$, at zero, up to order
$\zlma$; that is to say $\zgran a_{j}={{(2j-1)(2j-3)\zpu\zpu\zpu 1}\over {2^{j}\zpu j!}}\,$. Since the radius of
convergence of the power series $\zgran\zsu_{j=0}^{\zinf}a_{j}t^{j}$ is 1 and each $a_{j}>0$ , we have
$a_{0}=1\zmei \zf_{\zlma}\zmei (1-t)^{-1/2}$, $t\zpe [0,1)$, whence $(t-1)\zf_{\zlma}^{2}(t)+1\zmai 0$ and
$\zf_{\zlma}(t)\zmai 1$ when $t\zmai 0$ (both inequalities are obvious if $t\zmai 1$).

On the other hand if we set $\zf=t^{\zlma+1}R+\zf_{\zlma}$ then
$\zgran\zsu_{j=0}^{\zinf}t^{j}=(1-t)^{-1}=\zf^{2}=t^{\zlma+1}{\tilde R}+\zf_{\zlma}^{2}$ on $(-1,1)$; Therefore
$\zgran\zf_{\zlma}^{2}=t^{\zlma+1}S+\zsu_{j=0}^{\zlma}t^{j}$ where $S$ is a polynomial in $t$.
It follows, from that, the existence of a polynomial $\zl_{\zlma}$ of degree $\zlma$ such that
$(t-1)\zf_{\zlma}^{2}(t)+1=t^{\zlma+1}\zl_{\zlma}$.
\bigskip

{\bf Lemma 1.} {\it  For every $\zlma$ one has $(t-1)\zf^{2} +1= t^{\zlma +1}\zl_{\zlma}$ where $\zl_{\zlma}$ is
a polynomial of degree $\zlma$. Moreover $\zl_{\zlma}(t)\zmai 0$ and $\zf(t)>0$ for each $t\zpe {\Bbb R}$ if
$\zlma$ is even, and for any $t\zmai 0$ if $\zlma$ is odd.}
\bigskip

{\bf Proof.} It will suffice to show that $\zl_{\zlma}(t)\zmai 0$ and $\zf_{\zlma}(t)>0$ if $\zlma\zmai 2$ is even and
$t<0$. First we will prove, by induction on $\zlma$, the existence of a $\zd_{\zlma}>0$ such that $\zf_{\zlma}$ is
strictly decreasing on $(-\zinf, -1+\zd_{\zlma})$. Note that
$\zf_{\zlma}=a_{\zlma}t^{\zlma}+a_{\zlma-1}t^{\zlma-1}+\zf_{\zlma-2}=
a((2\zlma-1)t^{\zlma}+2\zlma t^{\zlma-1})+\zf_{\zlma-2}$ where $a>0$.

By induction hypothesis or because $\zf_{0}=1$, the polynomial $\zf_{\zlma-2}$ is decreasing on
$(-\zinf, -1+\zd_{\zlma-2})$, or on ${\Bbb R}$ if $\zlma=2$. But the derivative
$((2\zlma-1)t^{\zlma}+2\zlma t^{\zlma-1})'= ((2\zlma-1)\zlma t^{\zlma-1}+2\zlma(\zlma-1) t^{\zlma-2})<0$ on
$(-\zinf, -1]$, so $\zf_{\zlma}$ is strictly decreasing on some interval $(-\zinf, -1+\zd_{\zlma})$.

We show now that $\zf_{\zlma}(t)>(1-t)^{-1/2}>0$ if $t<0$. As $(1-t)^{-1/2}$ is strictly increasing, it is enough to
prove our result on $(-1,0)$. On this interval
$\zgran lim_{\zlma\zfl\zinf}\{ \zf_{\zlma}(t)\} =\zsu_{j=0}^{\zinf}a_{j}t^{j}=(1-t)^{-1/2}$.
But the series $\zgran\zsu_{j=0}^{\zinf}a_{j}t^{j}$ is alternating and the sequence
$\{a_{j}\zbv t\zbv^{j}\}_{j\zpe {\Bbb N}}$, whose limit is zero, strictly decreasing; then
$\zf_{\zlma}(t)>(1-t)^{-1/2}>0$ for $\zlma$ even.

Finally, if $\zf_{\zlma}(t)>(1-t)^{-1/2}>0$ for any $t<0$, a straightforward calculation shows that
$(t-1)\zf_{\zlma}^{2}(t)+1<0$, whence $\zl_{\zlma}(t)\zmai 0$ since $t^{\zlma+1}<0$. $\square$
\medskip

Recall that any polynomial $\zm$  which do not takes negative values has even degree and can be write
$\zm= \zm_{1}^{2} + \zm_{2}^{2}$, where $\zm_{1}$ and $\zm_{2}$ are polynomials of degree
$\zmei$ half of degree of $\zm$. Therefore by setting $k=\zlma +1$, $\za=\zf_{\zlma}$, $\zl_{\zlma}=\zm$,
$\zb_{1}=\zm_{1}$ and $\zb_{2}=\zm_{2}$ one has:
\bigskip

{\bf Corollary 1.} {\it For any odd natural number $k$ there exist three polynomials $\za,\zb_{1},\zb_{2}$, the
first one of degree $k-1$ and the other two with degree $\zgran\zmei {{k-1}\over 2}$, such that $\za(t)>0$
and $(1-t)\za^{2}(t) + t^{k}(\zb_{1}^{2}(t) + \zb_{2}^{2}(t))=1$ anywhere.}
\bigskip

Let us proof part (a) of theorem 1. Since topological degrees $\zmm 1$ may be represented by linear maps, we can
assume $k\zmai 1$. On ${\Bbb C}\zpor {\Bbb R}^{n-1} ={\Bbb R}^{n+1}$, endowed with coordinates
$(z,y)= (z, y_{1},..., y_{n-1})$ for which
$S^{n}=\{ (z,y);\, \zbv z\zbv^{2} +y_{1}^{2}+...+y_{n-1}^{2}=1\}$, we define
\medskip

\centerline{$F(z,y)=((\zb_{1}(\zbv z\zbv^{2}) +i \zb_{2}(\zbv z\zbv^{2}))z^{k}, \za(\zbv z\zbv^{2})y)$}
\medskip

\noindent where $\za$, $\zb_{1}$ and $\zb_{2}$ are as in corollary 1. Then $F(S^{n})\zco S^{n}$.

Set $S^{1}=\{ (z,0);\, \zbv z\zbv^{2}=1\}\zco S^{n}$. As $\za(t)>0$ for each $t\zpe {\Bbb R}$,
$F^{-1}(S^{1})=S^{1}$ and
$F$ preserves the orientation transversely to $S^{1}$. Hence the maps $F_{\zbv S^{1}}$ and
$F_{\zbv S^{n}}$ have the same topological degree, that is to say $k$.

By construction all the monomials of $F$ have odd degree $\zmei 2k-1$. Multiplying everyone of them by a suitable
power of $\zbv z\zbv^{2} +y_{1}^{2}+...+y_{n-1}^{2}$ the map $F$ becomes homogeneous of algebraic degree $2k-1$,
whereas $F_{\zbv S^{n}}$ do not change.

For proving (b), first we set ${\tilde\zl}_{\zlma}(t)=\zl_{\zlma}(t^{2})$ and
${\tilde\zf}_{\zlma}(t)=\zf_{\zlma}(t^{2})$. By lemma 1 we have
$(t^{2}-1){\tilde\zf}_{\zlma}^{2}(t)+1=t^{2\zlma+2}{\tilde\zl}_{\zlma}(t)$,
${\tilde\zf}_{\zlma}(t)>0$ and ${\tilde\zl}_{\zlma}(t)\zmai 0$ for any $t\zpe {\Bbb R}$. This allows us to find
out, for every natural number $\tilde k \zmai1$, three polynomials
${\tilde\za}, {\tilde\zb}_{1}, {\tilde\zb}_{2}$ such that ${\tilde\za}(t)>0$ and
$(1-t^{2}){\tilde\za}^{2}(t) + t^{2\tilde k}({\tilde \zb}_{1}^{2}(t) + {\tilde \zb}_{2}^{2}(t))=1$ anywhere.

Consider on ${\Bbb R}^{n+1}={\Bbb R}^{2r}\zpor {\Bbb R}^{n-2r+1}$
coordinates $(x,y)=(x_{1},...,x_{2r},y_{1},...,y_{n-2r+1})$. Let
$f:{\Bbb R}^{2r}\zfl {\Bbb R}^{2r}$ be a homogeneous polynomial
map of algebraic degree $2\tilde k$, sending $S^{2r-1}$ into
$S^{2r-1}$ with topological degree $k=\zmm2\tilde k$, which always
exists (see [3]) and $J:{\Bbb R}^{2r}\zfl {\Bbb R}^{2r}$ the
isomorphism given by $Jx=(-x_{2},x_{1},...,-x_{2r},x_{2r-1})$,
that is to say the canonical complex structure of ${\Bbb R}^{2r}$.
One defines (if ${\tilde k}=0$ just consider a constant map):
\medskip

\centerline{$F(z,y)=({\tilde\zb}_{1}(\zdbv x\zdbv^{2})f(x) + {\tilde\zb}_{2}(\zdbv x\zdbv^{2})Jf(x), {\tilde\za}
(\zdbv x\zdbv^{2})y)$}
\medskip

Then $F(S^{n}_{2r})\zco S^{n}$ and the same argument as in part (a), applied to
$S^{2r-1}=\{(x,0);\, \zdbv x\zdbv^{2} =1\}\zco S^{n}_{2r}$, shows that the topological degree of
$F:S^{n}_{2r}\zfl S^n$ equals $k$.
\bigskip

{\bf References}
\medskip

1. Baum, P.F.: Quadratics maps and stable homotopy groups of spheres, {\it Illinois J.Math.} {\bf 11} (1967),
586-595.

2. Golasi{\'n}ski, M. and G\'omez Ruiz. F.: Polynomial and regular maps into Grassmannians, {\it K-Theory} {\bf 26}
(2002), 51-68.

3. Wood, R.: Polynomial maps from spheres to spheres, {\it Invent. Math.} {\bf 5} (1968), 163-168.

4. Wood, R.: Polynomial maps of affine quadrics, {\it Bull. London Math. Soc.} {\bf 25} (1993), 491-497.

\end